\documentclass[a4paper,fleqn]{cas-sc}

\usepackage[authoryear,longnamesfirst]{natbib}
\usepackage{todonotes}
\usepackage{spverbatim}
\usepackage[acronym]{glossaries}
\usepackage{siunitx}
\usepackage{verbatimbox,caption,float,lipsum}
\usepackage{placeins}
\newfloat{Code}
\makeglossaries
\newacronym{apd}{APD}{action potential duration}
\newacronym{cdt}{CDT}{cardiac digital twin}
\newacronym{cep}{CEP}{cardiac electrophysiology}
\newacronym{erp}{ERP}{effective refractory period}
\newacronym{mi}{MI}{myocardial infarction}
\newacronym{uvc}{UVC}{Universal Ventricular Coordinates}
\newacronym{varp}{VARP}{Virtual Arrhythmia Risk Prediction}
\newacronym{vt}{VT}{Ventricular Tachycardia}
\newacronym{pcl}{PCL}{Pacyng Cycle Length}
\newacronym{ep}{EP}{ElectroPhysiology}
\newacronym{hpc}{HPC}{High Performance Computing}
\newacronym{rv}{RV}{Right Ventricle}
\newacronym{lat}{LAT}{Local Activation Time}
\newacronym{pp}{PP}{Pre-Pacing}
\newacronym{mt}{MT}{Maintenance}
\newacronym{ap}{AP}{Action Potential}

\newcommand{\aV}{\emph{auto-VARP}\xspace}
\newcommand{\cU}{\emph{carputils}\xspace}
\newcommand{\fC}{\emph{forCEPSS}\xspace}
\newcommand{\oC}{\emph{openCARP}\xspace}

\def\tsc#1{\csdef{#1}{\textsc{\lowercase{#1}}\xspace}}
\tsc{WGM}
\tsc{QE}
\tsc{EP}
\tsc{PMS}
\tsc{BEC}
\tsc{DE}

\begin{document}
\let\WriteBookmarks\relax
\def\floatpagepagefraction{1}
\def\textpagefraction{.001}


\shortauthors{P. Seghetti et~al.}

\title [mode = title]{AutoVARP -- A framework for automated reproducible inducibility testing in computational models of cardiac electrophysiology }                      



%

\author[1]{Paolo Seghetti}[
                        orcid=0000-0001-7511-2910]



\ead{paolo.seghetti@medunigraz.at}


\credit{Writing - Original draft preparation, Methodology, Software}


\author[1,2]{Matthias A. F. Gsell}[
orcid = 0000-0001-7742-8193
]
\credit{Resources, Data curation}
\author[1]{Anton Prassl}[
orcid = 0000-0001-7742-8193
]
\credit{Resources, Data curation}
\author[3]{Martin Bishop}[
orcid = 0000-0001-7742-8193
]
\credit{Conceptualization of this study, Methodology}


\author[1,2]{Gernot  Plank}[
orcid =0000-0002-7380-6908
   ]
\ead{gernot.plank@medunigraz.at}

\credit{Conceptualization of this study, Writing - Original draft preparation, Methodology, Software}

\affiliation[1]{organization={
    Gottfried Schatz Research Center, Division of Medical Physics and Biophysics, Medical  University of Graz},
    city={Graz},
    country={Austria}}
\affiliation[2]{organization={BioTechMed-Graz},
    city={Graz},
    country={Austria}}
\affiliation[3]{organization={
King's College London},
    city={London},
    country={United Kingdom}}




\begin{abstract}
\noindent\textbf{Background and Objective:} 
Simulations of \gls{cep} are gaining momentum beyond basic mechanistic studies,
as an approach for supporting clinical decision making.
The potential for \textit{in silico} technologies observed from the research community is immense, with studies demonstrating significantly improved therapeutical outcome with little to no additional burden for patients.
Studies replicating virtually induction protocols in post myocardial infarction patients are among the most reproduced and promise to identify non invasively ablation targets for therapeutical intervention.
Two main factors hinder the translation of these technologies from pure research to applications: virtually no reproducibility of results, and lack of standardized procedures. 
Inspired by a previously published virtual induction study by Arevalo et al. (2016), We address the issues of reproducibility and standardization providing \aV, a framework for standardization of virtual arrhythmia inducibility studies, built upon \oC and the \cU framework.\\
\textbf{Methods and Results:}
Standardization relies on the previously published forCEPSS framework and is ensured by defining the whole induction study with input files that can be easily shared to ensure reproducibility since the whole pipeline relies on open software.
Our approach also ensures numerical efficiency by separating the induction study into four stages: (i) pre-pacing with forCEPSS, (ii) S1 pacing tor each steady state, (iii) S2 induction with different extrastimuli, (iv) testing of sustenance of induced reentries.
We demonstrate the approach in a large virtual subject cohort to investigate numerical artifacts that may arise when improper setups are provided to perform virtual induction, and additionally showcase \aV in a biventricular mesh.\\
\textbf{Conclusions:}
\aV addresses effectively the current gap in standardization and reproducibility of results providing a uniform methodology that can be implemented even by non expert users. \aV is highly scalable and adaptable to markedly different anatomies. Although less flexible than \textit{in house} implementations it provides automated tools to share setups and does not require re-implementation of any process.  
\end{abstract}


\begin{highlights}
\item Open-source framework for standardizing inducibility studies in virtual hearts
\item High-level abstraction and workflow automation simplifies study execution 
and reduces user errors 
\item Efficient implementation of the induction protocol by staged execution
\item Facilitates verification of studies by reproduction due to reliance on freely available software \oC and \cU
\end{highlights}

\begin{keywords}
openCARP \sep
Standardized workflow \sep
Cardiac electrophysiology \sep
Ventricular tachycardia \sep
\end{keywords}

\maketitle

\section{Introduction}
Computational models of \gls{cep} offer a unique approach 
for gaining mechanistic insight into cardiac function in health and disease\cite{niederer2019:_models_cardiol}.
More recently, their scope has started to expand beyond basic research,
with applications geared towards medical device development~\cite{swenson20:_atp} in industry,
or even towards supporting clinical decision making~\cite{azzolin_personalized_2023,corral-aceroDigitalTwinEnable2020}.
An important class of clinical applications is the creation of patient-specific models,
often referred to as \glspl{cdt}\cite{baghirath24:_digital_twin},
for testing the inducibility of cardiac arrhythimas~\cite{sung22:_inFAT}. 
A specific point in case are \glspl{cdt} of patients surviving \gls{mi}
where the outcome of virtual inducibility tests is used to
stratify their risk of suffering from a \gls{vt}~\cite{arevalo16:_varp} 
or tailor \gls{vt} ablation therapies by guiding ablation targeting~\cite{prakosa18:_vaat}
in the real patient.

While conceptually appealing, the realization of such \emph{in silico} studies 
poses a number of major technical challenges.
First, the implementation of virtual inducibility testing studies 
based on protocols that mimic procedures used in the clinic 
for inducing a \gls{vt} in a patient constitutes a complex task.
Typically, this involves the definition of physiological parameters 
and their spatial heterogeneity, 
the selection of anatomical locations from where pacing stimuli are delivered, 
a stimulation protocol to define a limit cycle as an initial condition, 
a perturbation protocol for testing whether or not a \gls{vt} can be induced, 
and a suitable analysis over an observation period to determine 
whether the \gls{vt} is sustained, is similar to the clinical \gls{vt}, 
as well as the identification of \gls{vt} circuits and critical segments of the path
such as infarct isthmus exit sites~\cite{ciaccioVentricularTachycardiaSubstrate2025}. 
When executing such a study, a large number of simulations must be performed
and simulated outputs must be analyzed in a post-processing steps 
to automatically extract features and evaluate these against clinical criteria.
A sound implementation of such protocols, 
typically performed by individual specialized labs, 
is a demanding, time-consuming, non-standardized and error prone endeavor
that requires significant expertise beyond \gls{cep} in simulation engineering.
Owing to the large number of required simulations and the vast parameter space to sweep, 
a high level of processing automation is key to keep simulation workflows tractable
in terms of computational efficacy, robustness, and speed with minimal to no interactive processing.
Further, all processing steps and the used options must be encoded an audit trail 
to facilitate reproduction of a study. 
However, as the implementation of such protocols are often not shared nor made publicly available, 
scrutiny and reproduction is limited to the reporting laboratory.
As such, an independent external assessment of the findings of a study, whether for validation or falsification, remains unfeasible.

In this work we address this need for a standardized implementation of an \emph{in silico} inducibility testing protocol 
by providing an implementation of the \gls{varp} protocol~\cite{arevalo16:_varp}.
Our implementation of the \gls{varp} protocol, referred to as \aV, has been abstracted, 
highly automated, and generalized to be most broadly applicable and optimized for efficiency
to facilitate automated inducibility testing on large virtual cohorts, 
with an arbitrary number of pacing sites 
under flexible definitions of the stimulation protocol.  
The \aV pipeline is implemented in \href{https://www.python.org/}{Python} 
building upon the \cU framework~\cite{plank21:_opencarp},
using the standardized \fC tools for determining experimental limit cycles~\cite{gsell24:_forcepss}. 
The \aV simulation management is tailored for \oC, 
which is compatible with simulation tools \cite{vigmond08:_solvers} 
used in published studies using the \gls{varp} protocol. 
Post-processing and visualization builds upon the tools Meshtool~\cite{neic20:_meshtool} and Meshalyzer,
both included within the \oC framework.
The \aV implementation further supports the \oC bundle feature 
to simplify usage and reproduction.
The entire \aV implementation is freely available 
along with all simulation and processing tools,
thus ensuring a seamless reproduction of studies based on the \gls{varp} protocol.
The execution of a full blown inducibility study using \aV is showcased
by replicating a recent \gls{vt} simulation study 
which observed model-dependent stochastic effects with the \gls{varp} protocol~\cite{bishop25:_stochastic}.   

\section{Methods}
\begin{figure*}
    \centering
    \includegraphics[width=\textwidth]{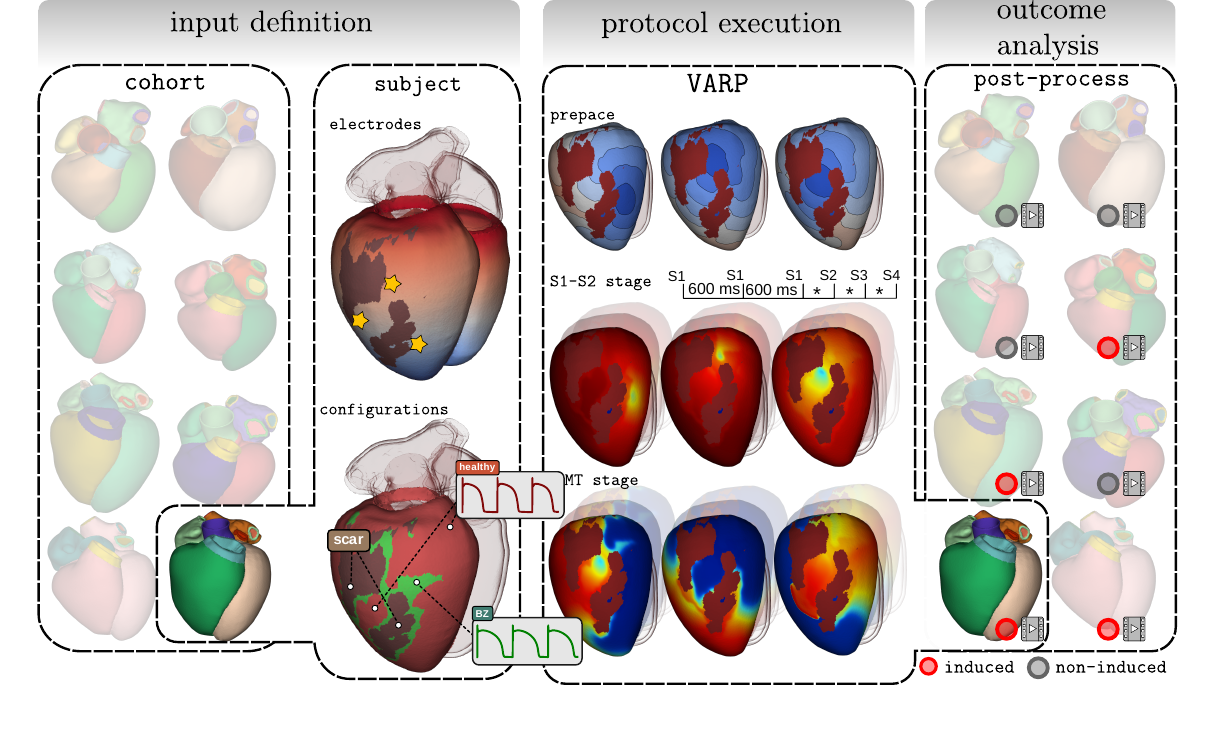}
    \caption{Outline of \aV workflow. A cohort of subject-specific anatomical meshes 
    with appropriate labeling is provided,
    along with definitions of subject-specific pacing locations 
    and electrophysiological properties characterizing normal healthy myocardium, 
    and impaired border zone tissue (left panels).
    The \gls{varp} inducibility test is executed for each subject in a staged workflow 
    comprising a pre-pacing (PP) stage to compute a limit cycle as a reference condition, 
    a S1 stage to stabilize the activation sequence to a prescribed \gls{pcl},
    a S2 stage to perform the inducibility test with variable coupling intervals,
    and a final maintenance (MT) stage to monitor reentrant activity and measure duration of \gls{vt} maintenance (middle panel).
    With the \gls{varp} protocol following an optimized pipeline (center).
    Simulation data of the maintenance stage are analyzed to generate tabular views on induction results 
    at the cohort level and for individual subjects and pacing locations. 
    All stages are stitched together to generate movies showing induction and maintenance 
    to facilitate lightweight review and documentation.}
    \label{fig:_overview}
\end{figure*}
An overview of the \emph{in silico} experimental workflow 
for assessing the vulnerability to \gls{vt} in post-\gls{mi} hearts
as implemented in \aV is given in Fig.~\ref{fig:_overview}.
A virtual cohort of cardiac models comprising tissue domains 
of heterogeneous \gls{cep} function 
along with a list of pacing sites and a defined pacing protocol 
are provided for inducibility testing.
The \aV pipeline processes these inputs to set up a sequence of simulations 
that apply the prescribed stimulation protocol at all pacing sites,
and monitor the electrical activity after the termination of the pacing protocol
to determine success or failure of \gls{vt} induction, 
and whether a \gls{vt} is maintained for a sufficiently long period according to predefined criteria.
For each simulation a binary output indicating success of induction is stored,
with the optional generation of movies for reviewing the entire induction experiment.

The implementation of the \aV pipeline builds upon open software frameworks 
integrated into the Python-based package \cU for managing the computational workflow.
The execution of \gls{cep} simulations relies upon \oC ~\cite{plank21:_opencarp} and \fC~\cite{gsell24:_forcepss},
\verb|meshalyzer| for visualization and \verb|ffmpeg| for movie generation. 

\subsection{Definition of Inputs}
The minimum set of inputs required for executing an inducibility test with the \aV pipeline includes 
a set of anatomical meshes, suitably tagged to specify domains of distinct \gls{cep} properties, 
a list of pacing sites for each anatomical mesh, 
a pre-defined pacing protocol to apply, 
and a dictionary encoding all experimental parameters in \verb|.json| format 
according to the \fC standard \cite{gsell24:_forcepss}. 

\subsubsection{Anatomical Meshes}
A cohort of anatomical meshes is provided in a format compatible with \oC~\cite{plank21:_opencarp}. 
Each subject of the virtual cohort must be stored in a subfolder within the virtual cohort directory, 
where the name of the subfolder corresponds with the subject identifier.
Each mesh must be appropriately labeled to delineate regions of distinct \gls{ep} properties
~(see Appendix \ref{sec:appendix:plan}). 

\subsubsection{Electrophysiological Properties}
\Gls{cep} tissue properties governing cellular dynamics or conductive properties 
are assigned to labeled regions within the anatomical mesh.
With \gls{varp} studies regions typically include: 
i) normal healthy myocardium, 
ii) functionally impaired border zone myocardium 
with altered repolarization and conduction properties\cite{arevalo16:_varp}, 
and iii) scarred tissue devoid of cardiomyocytes that is insulating and electrically non-excitable\cite{connollyComputationalRepresentationsMyocardial2016}.
However, our \aV implementation is flexible, able to incorporate an arbitrary number
of regions equipped with any \gls{ep} behavior that is supported by \oC. 
To ensure consistent behavior when using differently labeled regions,
a subject specific \verb|configurations| file can be provided in the folder containing the mesh 
to map each label to its respective \gls{ep} model (see Appendix \ref{sec:appendix:plan}).

\subsubsection{Pacing Locations}
An arbitrary number of pacing locations can be prescribed, 
either as Cartesian coordinates, as nodes of the anatomical mesh, or, 
most generically, based on anatomical coordinates such as \glspl{uvc}
which provide a single description of pacing sites 
that is equally applicable to all anatomical models of the cohort. 
This latter method is particularly useful when testing for arrhythmia inducibility in cohorts
of anatomically accurate virtual subjects, where Cartesian coordinates or mesh nodes cannot be used consistently.
Pacing locations are defined as electrodes in the \verb+electrodes+ section
of the input plan--file, or provided in a \verb|electrodes| file present in the subject directory~(see Appendix \ref{sec:appendix:plan}).

\subsubsection{Pacing Protocol}
All protocols defined in the planfile are executed for all subjects, to enforce standardization. 
Protocols may also be provided by an external \verb|protocols| file.
Protocols define the \glspl{pcl} used to reach tissue steady state and 
the pacing locations used, which may be different for each subject.
To ensure that each subject can use different pacing locations even when protocols are standardized,
each protocol refers to an electrode defined in the \verb+electrodes+ \fC section, 
which must be provided as a subject-specific \verb+electrodes+ file. As of now,
only one electrode can be defined per protocol since the standard S1-S2 induction
protocol always uses the same electrodes to deliver stimuli. 


\begin{figure}[h!]
        \centering
     \includegraphics[width=\textwidth]{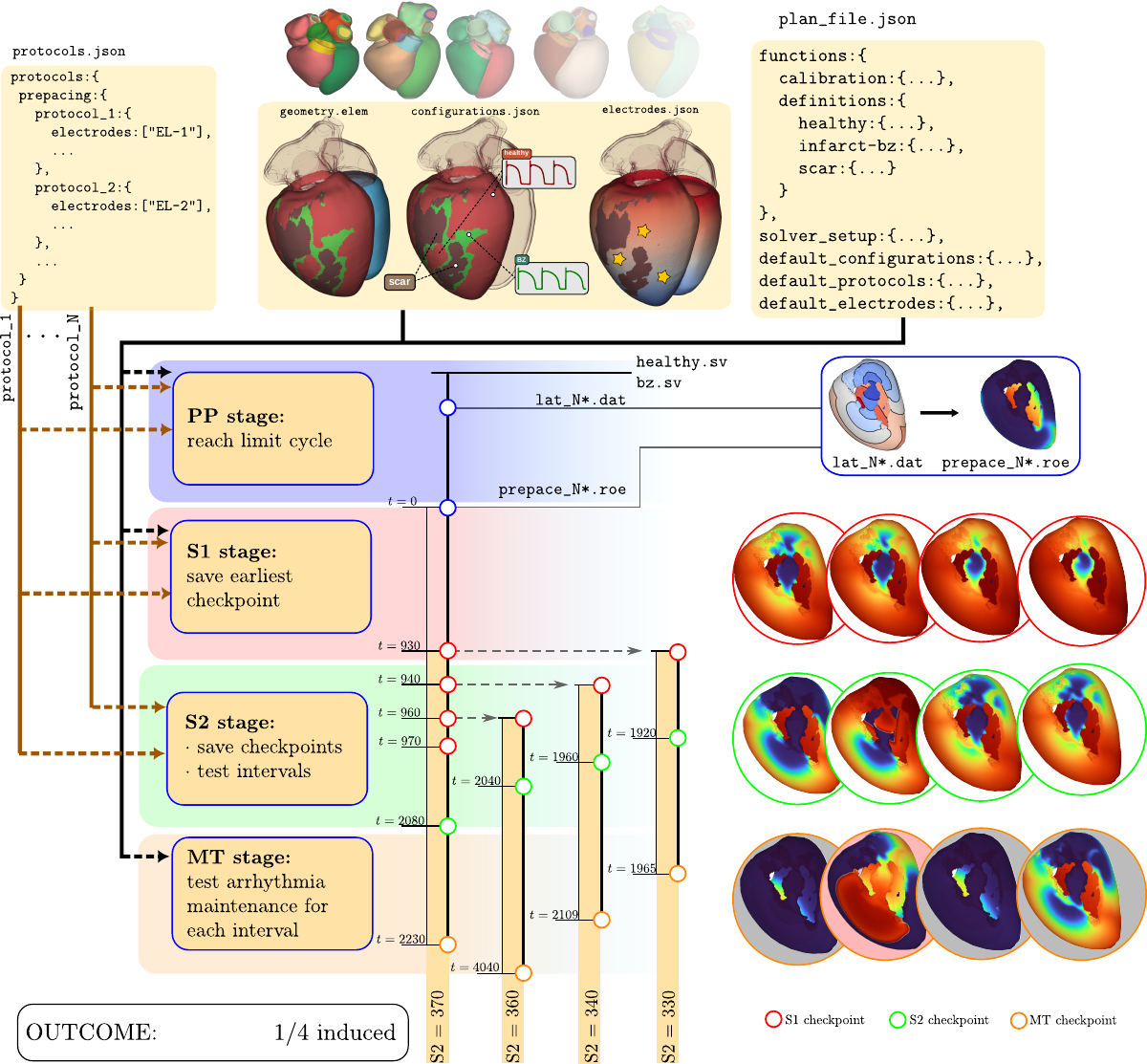}
    \caption{Complete \aV pipeline. Top: subject to be tested for induction, chosen in a cohort of other subjects. The subject specific \texttt{.json} files define electrodes and tissue configurations to be used in the induction study. The \gls{ep} settings are defined in the planfile, and the protocols to be tested are defined in the protocols file, which are general for all the virtual subjects in the cohort. 
    Bottom: each protocol from the protocols file is executed on each subject of the cohort. For each protocol, \aV stages are executed sequentially and checkpoint files are saved for each stage in order to restart simulations efficiently. 
    The figure shows snapshots of the membrane potential at each checkpoint instant, color coded for each stage. Each different S2 interval to be tested is independent of other S2 simulations and is restarted from a checkpoint generated using the largest S2 interval provided. At the end of each S2 simulations, a checkpoint file is saved and used to restart a MT simulation to test if the induction was successful. After the execution of \aV, the time stamps on MT checkpoints are used to establish if the induction was sustained or not.} 
    \label{fig:_stages}
\end{figure}

\subsection{Study definition and execution}

\aV automates the execution of a classic S1-S2 induction protocol in a cohort of virtual subjects using the \oC simulator.
The complex induction protocol is broken down into a sequence of simulation stages
to achieve computationally efficient execution and organizational tractability.
Stages are defined as simulations dedicated to a specific purpose 
(e.g.,\ achieving a tissue-level limit cycle as a reference condition 
approximating the normal sinus rhythm of a subject), 
aiming to avoid the repeated simulation of the same sequence. 
For instance, when testing different coupling intervals for delivering the S2 stimuli at a given pacing site, 
the entire preceding computation of limit cycle pacing and the delivery of the train of S1 stimuli 
is identical for all S2 stimuli, only the short period before delivering the first S2 stimulus will vary.
As such, storing a checkpoint at the shortest possible minimal S2 coupling interval 
to facilitate the continuation of a simulation safes compute time and vastly increases efficiency.  

Stages are executed sequentially for each subject in the cohort (Fig. \ref{fig:_stages}). 
The initialization phase follows the workflow implemented in the \fC framework 
to obtain efficiently a tissue/organ level limit cycle\cite{gsell24:_forcepss} (\gls{pp} stage),
followed by the S1 stage to entrain the tissue at the \gls{pcl} specified in the protocol, 
by the S2 stage to test inducibility by rapid pacing to generate unidirectional conduction blocks, 
and the final \gls{mt} stage to monitor reentrant activity and determine 
whether induction was successful and maintained for long enough according to pre-defined threshold criteria.
The organization in stages and the re-use of optimized protocols 
that exploit the specific capabilities of the \oC simulator yield a massive enhancement in efficiency.
For instance, alone the use of \fC to compute the limit cycle saves thousands of CPU hours
as compared to the brute force implementations used in previous studies~\cite{arevalo16:_varp,gsell24:_forcepss}.


In our implementation \aV sweeps sequentially all virtual subjects, 
simulates all the protocols provided in the input planfile, 
and executes all stages sequentially before moving to the next protocol (Fig. \ref{fig:_stages}). 
Before executing a simulation for a determined subject and stage,  
\aV checks if the simulation has already been performed by looking for the corresponding checkpoint file. 
If a checkpoint file is already present the simulation will be skipped, 
unless the \verb|--overwrite| flag is provided.
The checkpoint files generated by the S1 stage are saved in the \verb+checkpoints+ directory, 
associated with a \verb|.json| file containing metadata.
Due to the complex structure of the induction procedure, information must be retained between successive stages. 
As such, each call to \aV executing sequential stages must retain the inputs provided to the previous stage, to ensure
that all the necessary files to restart simulations are found. In the results section of this work, we report the command line used to generate each 
stage of the induction study, highlighting each new command that must be provided for each stage.

\subsubsection{Pre-Pacing stage}

The \gls{pp} stage computes a tissue/organ-scale limit cycle for each pacing site. 
The protocol is implemented in three different sub-stages according to the \fC pipeline, 
i.e.\ cellular limit cycle, an optional tissue tuning to calibrate conduction velocities, 
and determination of the limit cycle for a given electrode and \gls{pcl} \cite{gsell24:_forcepss}. 
The first sub-stage determines the \textit{cellular limit cycle} for each distinct model of cellular dynamics 
involved in the induction experiment, as defined in the \verb|functions| section 
(e.g.\ healthy myocardium and border zone tissue). 
In \aV, the number of pacing cycles and the \gls{pcl} are read from the plan file 
and are not directly an input. 
Since the same cellular dynamics models are shared across all subjects, 
the cellular limit cycle is computed only once and the corresponding state vectors 
are stored in \oC initial state files in a designated sub-directory as defined in the plan file.
The following optional sub-stage is a \emph{calibration of conduction velocities} by tuning conductivities 
to match prescribed orthotropic tissue level conduction velocities.
While this sub-stage was not considered in previous studies implementing \gls{varp}-like induction protocols \cite{arevalo16:_varp,prakosa18:_vaat,sung20:_personalized}, 
inclusion is recommended to avoid numerical artifacts leading to stochastic outcomes~\cite{bishop25:_stochastic}.
In reaction-diffusion simulations as implemented in \oC mesh resolution may affect conduction velocity 
due to numerical approximation errors \cite{gsell24:_forcepss}.
Specifically, this depends on the ratio between conduction velocity and mesh resolution
and is particularly critical for slow conduction.
For small ratios associated with slow conduction wave fronts are spatially undersampled,
leading to extremely non-physiologically short wave lengths and artificial conduction blocks,
rendering an induction experiment unreliable. 
By providing the \verb+--tissue-tuning+ flag, each \gls{cep} function definition 
is calibrated by tuning conductivities until wave front conduction velocity matches the values prescribed in the plan file.
In absence of this flag, conductivity settings will be used as provided, 
and associated conduction velocities will be measured and documented in the plan file.
Once conductivities are determined, 
\textit{tissue-level limit cycle} is computed by executing the \verb|limit_cycle| program of \fC,
using one of the prepacing options offered within \fC, 
exposed through the \verb+--gen-lat+ and \verb+--lim-cyc+ options. 
These approximate the tissue/organ stable limit cycle for a given electrode and \gls{pcl} 
based on a cellular limit cycle by accounting for the activation sequence and diffusion effects 
with variable degrees of computational cost and accuracy. 
The global electrophysiological state of the tissue after this procedure
approximates a tissue state corresponding to pacing the preparation at the tissue/organ model for a large number of cycles.
For each subject and each pacing site, the tissue/organ limit cycle is stored as an \oC  checkpoint (\verb+.roe+ files), 
automatically labeled with a specific ID to facilitate automated checkpoint selection at the following stages of \aV 
(see Appendix \ref{sec:appendix:checkpoints}). 
Checkpoints and pacing-induced activation maps are saved in the \verb+checkpoints+ folder 
located inside each subject directory, to be used then to initialize the following S1 stage. 
In addition to initializing the computational experiments, the \gls{pp} stage
can also  be used to generate restitution curves for the \gls{ep} tissues used in the simulations, 
by providing the flag \verb|--plot-restitution|. 
In this case, \aV shows the restitution curves for all the tissues involved in the experiment and the respective \glspl{apd}, and saves the plots as images, using \fC. 
This feature is useful to determine a priori the intervals S2 to be tested for induction. 
Indeed, infarct-related anatomical reentry is triggered by unidirectional conduction block caused by
dispersion of repolarization between healthy tissue and infarct border zone. By visually inspecting the 
restitution curves, the user can determine the extrastimulus intervals that will excite healthy tissue
but will be blocked by the border zone, still refractory.

\subsubsection{S1 stage}
The S1 stage is performed to entrain the preparation by pacing at the prescribed S1 \gls{pcl}, 
and to ascertain that a true stable limit cycle has been achieved, i.e.\ at any point in space 
cellular dynamics traverses the state space following the same trajectory.
While up to nine S1 pacing stimuli can be delivered, 
the actual choice of cycles depends on the chosen mode of pre-pacing.
For \gls{lat}-based \gls{pp} modes the number of stimuli required to arrive at a S1 limit cycle is less \cite{gsell24:_forcepss}.
Typically, as the same \gls{pcl} is used at \gls{pp} and S1 stage, a single S1 cycle suffices. 
To save computational costs, the S1 stage terminates by checkpointing the tissue/organ state,
starting at the shortest S2 coupling interval to be considered. Since during both S1 and S2 stages the protocol is enforced (S1-S2 protocol), electrodes are used to define the name of the checkpoint files, in addition to stage and \gls{pcl} (see appendix \ref{sec:appendix:checkpoints}). 

Typically, the effective refractory period of the tissue at the stimulus location is chosen,
and provided to \aV as the smallest coupling interval to test, with the input argument \verb|--CI-array|. 
To cover all possible S2 coupling intervals, the argument \verb|--CI-array| is provided as a comma separated user-defined list of pacing intervals to test in the S2 stage.
As long as the same minimum S2 in \verb+--CI-array+ is used between different executions of the program, 
the S1 stage will not be recomputed and \aV will automatically restart from the S2 stage. 
Note that if testing for S2 intervals shorter than what previously defined in \verb+--CI_array+ is desired, 
the S1 stage must be recomputed. 

\subsubsection{S2 stage}

The S2 stage performs the actual inducibility test. 
Whether an arrhythmia such as a \gls{vt} is inducible or not depends on
the topology of the substrate and the electrophysiological properties,
the location of the pacing site relative to the substrate, 
and, most importantly, the timing of stimulus delivery.
Induction according to the \gls{varp} protocol mechanistically builds on a creating an unidirectional block,
generated by the restitution properties and dispersion of repolarization in the substrate\cite{enriquezPathophysiologyVentricularTachyarrhythmias2017}. 
With rapid stimulation such a block is provoked at the interface between tissues of different \glspl{erp},
as between healthy and impaired border zone tissue, if a S2 stimulus is delivered prematurely, that is, S2$<$S1.
The window for prematurity is bounded by the \glspl{erp} of healthy tissue, $\text{ERP}_{\rm h}$, and impaired tissue, $\text{ERP}_{\rm{bz}}$.
The shorter effective refractory period within healthy tissue 
is the lower bound of the window 
as S2$<\text{ERP}_{\rm{h}}$ would lead to loss of capture.
The upper bound is due to $\text{ERP}_{\rm{bz}}$ 
as with S2$>\text{ERP}_{\rm{bz}}$ no unidirectional block is formed anymore.
These bounds can be estimated based on the \gls{apd} restitution curves computed within \fC (e.g., see Fig. \ref{fig:prepace}).
For the electrophysiological settings used within \gls{varp} as published\cite{arevalo16:_varp}
the prematurity window falls within \SI{280}{\milli \second} and \SI{350}{\milli \second}.
The upper bound is an estimate only as the effect of conduction velocity restitution is neglected.
That is, for instance, S2$<$\SI{350}{\milli \second} could fail to create a unidirectional block at an interface
between healthy and border zone tissue due to conduction slowing in healthy tissue, 
postponing the arrival time at the interface, providing a longer time window for border zone tissue to recover.

Therefore, in general, it is necessary to test various S2 intervals 
to gauge the arrhythmogenic potential of a determined substrate. 
There are various S2 intervals able to induce a \gls{vt}, and, 
depending on a chosen S2, the induced \glspl{vt} may behave differently, or, 
in topologically more complex substrates, \glspl{vt} traversing topologically different circuits may be elicited.
This implies that several S2 stimuli need to be tested, each continuing from the same S1 simulation.
To save computational resources, each simulation associated with a different extra-stimulus interval 
is restarted from the same S1 checkpoint saved in the S1 stage.
If more than one S2 intervals are provided for testing (e.g., \SI{280}{\milli \second},\SI{300}{\milli \second},\SI{350}{\milli \second}), 
the S2 simulation with the largest S2 interval is used to generate starting checkpoints for all the remaining smaller S2 intervals, saving further simulation time (refer to Fig.\ref{fig:_stages}). 
During this S2 simulation, all other starting S1 checkpoint files for the shorter S2 intervals are saved, since no stimulus is being delivered. After the checkpoints are saved, the remaining S2 simulations are started.
Both the number of S2 (one to nine) and S2 intervals can be adjusted by the user, 
by using the arguments \verb+--S2-cycles+ and \verb+--CI-array+.
Additionally, it is possible to deliver a decremental train of stimuli 
by providing a decrement interval with the argument \verb|--decrement-S2|. 
Each S2 simulation is carried out until the end of the S2 stimulation protocol, 
plus an additional waiting time lasting for the duration of a single S2 interval. 
At the end of each S2 simulation, a checkpoint is generated and saved in the \verb|checkpoints| directory
to be used as a restarting checkpoint for continuing the simulation at the \gls{mt} stage of the \aV program.
\subsubsection{Maintenance stage}
The final simulation stage performed by \aV is used to analyze the induced activity
with regard to success or failure of induction,
and -- for successful induction -- for how long the arrhythmia is sustained.
Each S2 simulation is continued by restarting from the corresponding checkpoint. 
The evolution of \gls{cep} activity in the tissue/organ preparation is monitored 
over a maximum observation period of \SI{2}{\second}, 
or until no active wave front propagation is detected over a \SI{150}{\milli \second} time window.
To determine wether or not \gls{cep} activity in the preparation has ceased, \aV uses the \verb|sentinel| option of 
\oC to detect \gls{ap} upstroke. 
This ensures that unsuccessful induction tests are stopped without wasting computational resources.

Depending on the criterion used for determining maintenance, 
the lengths of the observation periods can be modified by the user by providing the argument \verb|--MT-duration|.
The default of of two seconds has been chosen according to previous studies\cite{arevalo16:_varp}
Regardless of the duration of the simulation, a labeled checkpoint is saved before quitting
where the label indicates the simulated MT duration (see appendix \ref{sec:appendix:checkpoints}), used then in the post-processing phase 
to evaluate the outcome. In the post processing phase, the MT checkpoint file name is used to determine the electrode and extrastimulus intervals that are associated with the duration of the \gls{mt} simulation, which is used to determine if the stimulation caused reentrant arrhythmia.

\subsubsection{Post processing options}
\label{sec:pp_opts}
When evaluating the outcome of inducibility experiments in larger cohorts, 
owing to the large number of simulations to review a manual inspection of outcomes 
is time consuming and error prone. 
Therefore, \aV is provided with two post processing options to facilitate automated analysis. 

The option\hspace{0.5cm}\verb+--status+ informs on the state of the \aV pipeline for each subject 
found in the cohort directory. 
For each subject, the status is determined by examination of checkpoints
and printed as table where rows display the protocols used in the pipeline 
and the columns display the progress in terms of completed stages. 
For completed MT stages the table also displays 
how many sustained \glspl{vt} were induced for each electrode, 
by using a time threshold, provided with the argument \verb|--mt-threshold|,
to discriminate sustained and non-sustained \gls{vt}. 
Tables are saved in comma-separated format in the \verb+status_tables+ folder.

The option \hspace{0.5cm}\verb+--movies+ generates animations of the spatio-temporal evolution 
of transmembrane voltages using Meshalyzer\cite{plank21:_opencarp} for each S2 interval tested. 
so that if data are removed or not found the command fails. 
Meshalyzer uses a user-provided \emph{.mshz} file in the root directory of the script 
to load specific configurations of the scene used for rendering the movies. 
\aV is provided with a default settings file, that can be edited using Meshalyzer.

The option \hspace{0.5cm}\verb|--bundle| is used to generate a compressed folder 
containing all the files necessary to reproduce the \gls{varp} experiment. 
This feature previously implemented in \oC \cite{houillon_facilitating_2024} 
facilitates replication of simulation results by generating a self-contained,
an automated upload to the \href{https://radar.kit.edu/radar/de/home}{RADAR} server, 
or a manual provision of the bundle in any open repository (e.g., \href{https://zenodo.org/}{Zenodo}).

\section{Results}
We showcase the execution of \aV for testing the inducibility of \glspl{vt},
as previously used in studies investigating arrhythmogenesis~\cite{arevalo16:_varp,sung22:_inFAT},
in two scenarios of increasing complexity: 
i) a virtual cohort comprising 36 subjects modeled as 3D tissue slabs with simple scar geometries, and, 
ii) a cohort of a single biventricular setup with \gls{uvc} parametrization. 
The simple slab cohort with scars traversed by an isthmus were generated 
to replicate a recent study or ours~\cite{bishop25:_stochastic} 
that showed the generation of numerical artifacts due to inappropriate mesh resolution,
leading to stochastic outcomes of inducibility tests.
In the geometrically simplified models the overall geometric shape of scar and isthmus were kept fixed,
only details of the implementation were altered. 
We chose the width of infarct border zone as \SI{1}{\milli \meter} or \SI{2}{\milli \meter}, 
discretized using 6 different mesh resolutions, and imposed 3 fiber fields 
differing in their angle relative to the isthmus, 
yielding a large cohort of subjects suitable for showcasing an entire \aV workflow including  post-processing 
(see Fig.~\ref{fig:results_fig} and Tab.~\ref{tab:status} for details).
The more complex biventricular setups demonstrate that the \aV workflow is sufficiently abstracted 
to be universally applicable to any cohort of anatomical meshes, regardless of their shape or labeling.

\subsection{Setting up a new \gls{varp} experiment}
Initiating an \aV experiment requires: 
i) a planfile \verb|planfile.json|,
ii) a \verb+cohort+ folder containing all virtual subjects (36 isthmus geometries), 
and, 
iii) a protocols file \verb|varp_protocols.json|. 
The planfile contains the definition of the three \gls{ep} functions 
that define healthy, border zone and scar tissue, their standard configurations, and default solver settings.
Each subject subfolder contains the geometry files (\verb+.elem+, \verb+.pts+, \verb+.lon+), 
and, optionally, 
experiment files (\verb+.json+ files containing electrodes and configurations) to customize protocols for individual subjects
by overruling default settings in the planfile, 
and an anatomical reference frame in the form of e.g.\ universal ventricular coordinates file (\verb+.uvc+) 
in cases where electrode locations are defined in relative anatomical terms
to ensure electrode locations are place equally throughout a cohort.
For the simple 3D slab cohort where the geometrical shape and anatomical labeling of all virtual subjects is the same
the provision of subject specific configurations and electrodes files is not required, 
as electrodes and linking of anatomical labels to electrophysiological function used in the induction study 
can be uniquely defined throughout the entire cohort in the planfile. 
Electrodes are geometrically defined in Cartesian coordinates and arranged around the isthmus, shown in Fig.\ref{fig:results_fig}a.

The command line options pertaining to the experimental setup are:
\begin{enumerate}[{}]
    \item \hspace{0.5cm}\verb+--plan+ to provide the path to the plan file defining solver settings, 
    \gls{cep} function settings, and default protocols, electrodes, and configurations.
    \item \hspace{0.5cm}\verb+--cohort-dir+ or \verb+--case-ID+ to provide the path to the cohort directory containing the subjects folders, 
    or the path to the subject folder to restrict processing to one selected subject.
    \item \hspace{0.5cm}\verb+--protocols+ to provide the path to the protocols file 
    defining all the protocols to be executed on the cohort. 
    \item \hspace{0.5cm}\verb+--configurations+ to provide a configurations file 
    in each subject subfolder (name must be the same for each subject) 
    linking anatomical domain to electrophysiological function. 
    \item \hspace{0.5cm}\verb+--electrodes+ to provide a subject-specific electrodes file 
    (name must be the same for each subject). 
\end{enumerate}
The options \verb+--plan+ and \verb+--cohort-dir+ are mandatory, 
providing the rules for executing the induction experiments, 
and the cohort of subjects to operate on.
All other options are optional, if not provided default options included in the plan file will be used.
As such, the entire \gls{varp} experiment could be executed automatically by executing a single command:\vspace{3mm}\\
\verb+python auto-varp.py --plan planfile.json --cohort-dir cohort+\vspace{3mm}\\
In general though, \aV will be executed stage-by-stage, including a monitoring of intermediate results.
This is of particular importance for actual clinical studies 
where a number of patient-specific anatomical meshes will be processed.
In an unsupervised execution of \aV errors in the provided input files may remain undetected, 
which may entail wasting a vast amount of computational resources.

For stage-wise execution the provision of \aV stage-specific arguments is incremental, 
requiring the specification of all arguments used at preceding stages, 
to ensure a correct continuation of the protocol from previous checkpoints and, thus,
the definition of a unique experimental protocol through the entire tree of possible experiments.  
As such, the stage-by-stage execution is demonstrated in the following for the simple slab cohort, 
using the full spectrum of available options.

\subsection{Pre-Pacing stage - establishing a stable limit cycle}
\label{res:pp}
We execute the \gls{pp} stage using the most complete pre-pacing option, 
i.e.\ \gls{lat}-based distribution of initial states followed by a single simulated cycle 
using the \oC reaction-diffusion monodomain model.
This is triggered by using the options \verb|--gen-lat ek --lim-cyc lat-1|. 
Since we are interested in studying the effect of mesh resolution in the induction protocol, 
we do not perform tuning of conductivities. 
As all 36 meshes in the slab cohort are tagged using the same labels, 
\verb|1| for healthy tissue, \verb|2| and \verb|3| for border zone tissue, and \verb|4| for scar tissue, 
the provision of a configurations file is not required when launching \aV. 
The \gls{pp} stage is executed for the entire cohort of subjects by:\vspace{3mm}
\\
\small
    \verb|python auto-varp.py --plan planfile.json --cohort-dir cohort --protocols varp_protocols.json \|\\
    \verb|                    --stage PP --gen-lat ek --lim-cyc lat-1| 
\normalsize
\vspace{3mm}

\noindent After execution, in each subject folder, a \verb|checkpoints| folder is created 
containing the activation sequence file (\verb+.dat+) generated by the \verb|--gen-lat| option 
and the checkpoint file (\verb+.roe+) generated by the \verb|--lim-cyc| option, for each protocol provided. 
In this case, eight activation sequence files and eight checkpoint files are generated, 
each named with a unique ID generated including the stage, protocol used, subject, and employed prepacing-mode 
to distinguish between different simulations (see appendix \ref{sec:appendix:checkpoints}).  
Additionally, the cellular limit cycles of ionic models defined in the planfile 
have been saved as initial state vectors in the \verb|init| folder. 
In this example, two initialization states are created corresponding to the healthy and border zone tissues.

\begin{figure}
    \centering
    \includegraphics[width=\linewidth]{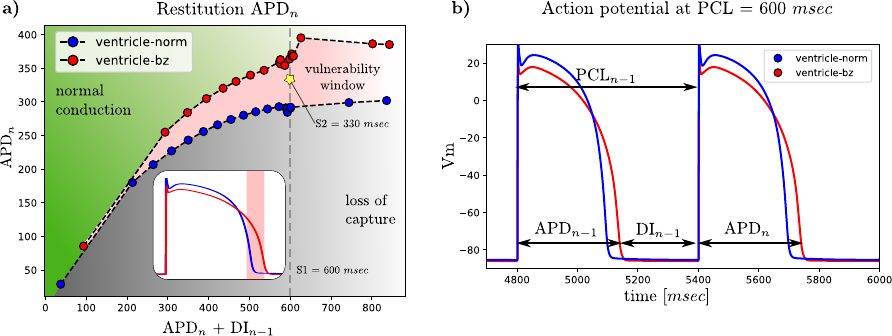}
    \caption{Single cell behavior of the ionic models considered in the \aV induction study. 
    a) At the S1 \gls{pcl} a marked difference in \gls{apd} between healthy and border zone tissue is noted, 
    opening a window of vulnerability. 
    By delivering a premature stimulus within this window, 
    the healthy tissue will be excited but not the border zone tissue, causing a unidirectional conduction block. 
    b) Last two \glspl{ap} obtained with a \gls{pcl} of \SI{600}{\milli \second} for 100 beats. 
    At the cellular limit cycle $PCL_{n} = APD_n + DI_{n}$ holds.  }
    \label{fig:prepace}
\end{figure}

By using the \verb|--plot-restitution| option, the restitution curves of healthy and border zone tissue are visualized and  used for determining a window of interest for the induction of unidirectional conduction block (Fig.\ref{fig:prepace}). 
For the chosen S1 \gls{pcl} of \SI{600}{\milli \second}, the window of vulnerability starts at  \SI{280}{\milli \second} and ends at \SI{350}{\milli \second}. 
For S2 coupling intervals lower than \SI{280}{\milli \second}, both tissues are refractory and loss of capture occurs, 
whereas for values larger than \SI{350}{\milli \second} both tissues are excitable and support conduction.
Thus, at the S2 stage a coupling interval of \SI{330}{\milli \second} is chosen to induce a conduction block at the interface between normal and border zone tissue.

\subsection{S1 stage -- preconditioning for inducibility testing}
Up to 9 consecutive S1 pacing stimuli can be delivered using the input argument \verb|--S1-cycles|
to entrain the tissue at a given S1 \gls{pcl}.
As \gls{pp} and S1 stage are using the same cycle length, 
and the \gls{pp} used the most accurate method including the simulation of a full cycle 
with a reaction-diffusion monodomain model, 
all setups will be at the S1 limit cycle, or very close to it. 
As such, executing a single S1 cycle suffices to ensure that the S1 limit cycle has been reached. 
To perform the S1 stage the following command is issued:\vspace{3mm}\\
\small
    {\color{gray}\verb|python auto-varp.py --plan planfile.json --cohort-dir cohort --protocols varp_protocols.json\|}\\
    \verb|    --stage S1|{\color{gray}\verb| --gen-lat ek --lim-cyc lat-1\|}\\
    \verb|    --S1-cycles 1 --CI-array 280,330|
\normalsize
\vspace{3mm}

Note that command line parameters used at the preceding \gls{pp} stage (light gray) are repeated 
to ascertain continuation of the experiment using the checkpoint corresponding to the chosen experimental history. 
The additional parameter \verb|--CI-array 280,330| is given 
for minimizing computational expenses at the following S2 stage. 
The shortest coupling interval to be tested at the S2 stage is determined as the minimum value of the \verb|--CI-array|,
which is chosen based on the \gls{apd} witnessed in healthy tissue at the S1 \gls{pcl}.
Shorter coupling intervals may be feasible if a shorter S1 \gls{pcl} is selected, yielding a shorter \gls{apd}.
By providing \SI{280}{\milli \second} as an additional interval, we ensure that the minimum capture interval is saved as a checkpoint and can be used as a starting checkpoint in future induction studies. 


\noindent After execution, 16 S1 checkpoint files are generated, two for each protocol, at the times specified by \verb|--CI-array|, 
and are saved in the \verb|checkpoints| folder of each virtual subject. 
Additionally, the output of each simulation including spatio-temporal transmembrane voltage data and log files 
are saved in the \verb|sim_outputs| folder located in the root directory of the program.
Representative snapshots of transmembrane voltage distribution at the start and at the end of the S1 stage are shown in Fig.~\ref{fig:results_fig}b for a single subject of the cohort.
Note that the S1 stage is stopped at the minimum coupling interval provided (\SI{280}{\milli \second} in this example).

\subsection{S2 stage - premature stimulus delivery}

At the S2 stage only the number of stimuli to be delivered with the given S2 coupling interval is provided.
Based on data measured at the \gls{pp} stage (Fig.~\ref{fig:prepace}),
using a S2 \gls{pcl} of \SI{330}{\milli \second} will induce a conduction block 
at the interface between healthy and border zone tissue.
This unidirectional block prevents the wave front to directly propagate into the proximal entrance of the isthmus and, thus, sets up the condition for a \gls{vt} circuit to form 
as the wave front propagates through the outer loop around the scar, 
to reenter the isthmus then later through the distal entrance of the isthmus.
Hence, we only provide the option \verb|--S2-cycles 1| to probe inducibility of the \gls{vt} circuit
from all pacing sites
\vspace{3mm}\\
\small
    {\color{gray}\verb|python auto-varp.py --plan planfile.json --cohort-dir cohort --protocols varp_protocols.json\|}\\
    \verb|    --stage S2|{\color{gray}\verb| --gen-lat ek --lim-cyc lat-1\|}\\
    {\color{gray}\verb|    --S1-cycles 1 --CI-array 280,330\|}\\
    \verb|    --S2-cycles 1|
\normalsize
\vspace{3mm}


\noindent For the sake of computational efficiency, 
the simulation corresponding to the longest S2 coupling interval provided in \verb|--CI-array|
generates the checkpoints for all shorter S2 intervals to be tested (see Fig.~\ref{fig:_stages}).
The S2 stage simulates activity in response to the S2 stimulus for a full S2 \gls{pcl}, 
and terminates then by storing checkpoints in the \verb|checkpoints| directory of each subject.
Representative snapshots of transmembrane voltage distribution 
at the beginning and at the end of the S2 stage are shown in Fig. \ref{fig:results_fig}b for a single subject of the cohort.

\subsection{MT stage -- testing arrhythmia maintenance}

No further input options are needed at the MT stage when evaluating the maintenance criterion, 
typically requiring reentrant activity to be self-sustained for $>$\SI{2}{\second}.
Simulations are re-started from the terminal S2 checkpoints, without delivering any further stimuli,
and run over the full prescribed time window of \SI{2}{\second}
by executing:\vspace{3mm}\\
\small
    {\color{gray}\verb|python auto-varp.py --plan planfile.json --cohort-dir cohort --protocols varp_protocols.json\|}\\
    \verb|    --stage MT|{\color{gray}\verb| --gen-lat ek --lim-cyc lat-1\|}\\
    {\color{gray}\verb|    --S1-cycles 1 --CI-array 280,330\|}\\
    {\color{gray}\verb|    --S2-cycles 1|}
\normalsize
\vspace{3mm}

\noindent At the MT stage activity during simulation is constantly monitored.
If no wave front propagation is detected over a given observation period of \SI{150}{\milli \second} 
simulations terminate prematurely, storing a checkpoint.
All checkpoints are labeled with a time stamp indicating the instant of termination
to facilitate the evaluation in a post-processing step. 
Prematurely terminated simulations are classified as non-inducible or non-sustained,
and completed simulations as maintained.
Figure~\ref{fig:results_fig}b shows the outcome of \aV for subject \verb|1mmbz.300um.f90| of the virtual cohort. Two induction sites resulted in sustained reentry, 
whereas all other simulations terminated prematurely as electrical quiescence, 
i.e.\ absence of propagating wave fronts,
was detected over an observation period of \SI{150}{\milli \second}.

\subsection{Post-processing stage - progress monitoring and outcome evaluation}

Executing \aV on a cohort of 36 subjects where in each subject induction was tested 
from 8 electrode locations with a single S2 coupling interval using four processing stages 
produced a large number of outputs.
In this study the MT stage alone produced 288 simulations.
Outputs must be analyzed to monitor the progress status of the \aV workflow,
and for completed experiments the outcome must be evaluated.
Considering the vast number of simulations, manual inspection is prohibitive.
Thus, \aV provides a \verb|--status| option to print workflow progress for each subject in a tabular format. 
Progress status can be interrogated then by issuing
\vspace{3mm}\\
\begin{minipage}{\textwidth}
\small
    {\color{gray}\verb|python auto-varp.py --plan planfile.json --cohort-dir cohort --protocols varp_protocols.json\|}\\
    \verb|    --stage MT|{\color{gray}\verb| --gen-lat ek --lim-cyc lat-1\|}\\
    {\color{gray}\verb|    --S1-cycles 1 --CI-array 280,330\|}\\
    {\color{gray}\verb|    --S2-cycles 1\|}\\
    \verb|    --status --mt-thr 1000|
\normalsize
\end{minipage}
\vspace{3mm}

\noindent where the argument \verb|--status| triggers the status output, 
and \verb|--mt-thr| is provided to define the cutoff value in \si{\milli \second}
for classifying the outcome of the inducibility protocol as 
non-inducible, inducible but not sustained and sustained.
An exemplar tabular output for testing inducibility in subject \verb|1mmbz.300um.f90| of the cohort
is shown in Tab.~\ref{tab:status}. 
Each row of the table corresponds to one tested electrode and displays the checkpoints computed on each column.

\begin{table}
    \centering
    \resizebox{\linewidth}{!}{
    \begin{tabular}{cccccc}
    \hline
        protocol and electrode & gen-lat & lim-cycle & S1 stage(BCL-checkpoint) & S2 stage (BCL-S2-checkpoint) & MT stage (BCL-S2-duration)\\
        \hline
         prepace-S1200-RAD,S1200-RAD:&  gen-lat&  lim\_cyc-0600&  S1-0600-280/S1-0600-330&  S2-0600-0330-660& MT-0600-0330-02000/over 1000 msec: 1/1\\
         prepace-S0130-RAD,S0130-RAD:&  gen-lat&  lim\_cyc-0600&  S1-0600-280/S1-0600-330&  S2-0600-0330-660& MT-0600-0330-00149/over 1000 msec: 0/1\\
         prepace-S0300-RAD,S0300-RAD:&  gen-lat&  lim\_cyc-0600&  S1-0600-280/S1-0600-330&  S2-0600-0330-660& MT-0600-0330-00149/over 1000 msec: 0/1\\
         prepace-S0430-RAD,S0430-RAD:&  gen-lat&  lim\_cyc-0600&  S1-0600-280/S1-0600-330&  S2-0600-0330-660& MT-0600-0330-00150/over 1000 msec: 0/1\\
         prepace-S0600-RAD,S0600-RAD:&  gen-lat&  lim\_cyc-0600&  S1-0600-280/S1-0600-330&  S2-0600-0330-660& MT-0600-0330-02000/over 1000 msec: 1/1\\
         prepace-S0730-RAD,S0730-RAD:&  gen-lat&  lim\_cyc-0600&  S1-0600-280/S1-0600-330&  S2-0600-0330-660& MT-0600-0330-00149/over 1000 msec: 0/1\\
         prepace-S0900-RAD,S0900-RAD:&  gen-lat&  lim\_cyc-0600&  S1-0600-280/S1-0600-330&  S2-0600-0330-660& MT-0600-0330-00149/over 1000 msec: 0/1\\
         prepace-S1030-RAD,S1030-RAD:&  gen-lat&  lim\_cyc-0600&  S1-0600-280/S1-0600-330&  S2-0600-0330-660& MT-0600-0330-00150/over 1000 msec: 0/1\\
    \end{tabular}
    }
    \caption{Tabular view of workflow progress status generated from \aV for subject \texttt{1mmbz.300um.f90}
    using \texttt{--status}. } 
    \label{tab:status}
\end{table}

\begin{figure*}[h!]
    \centering
    \includegraphics[width=\textwidth]{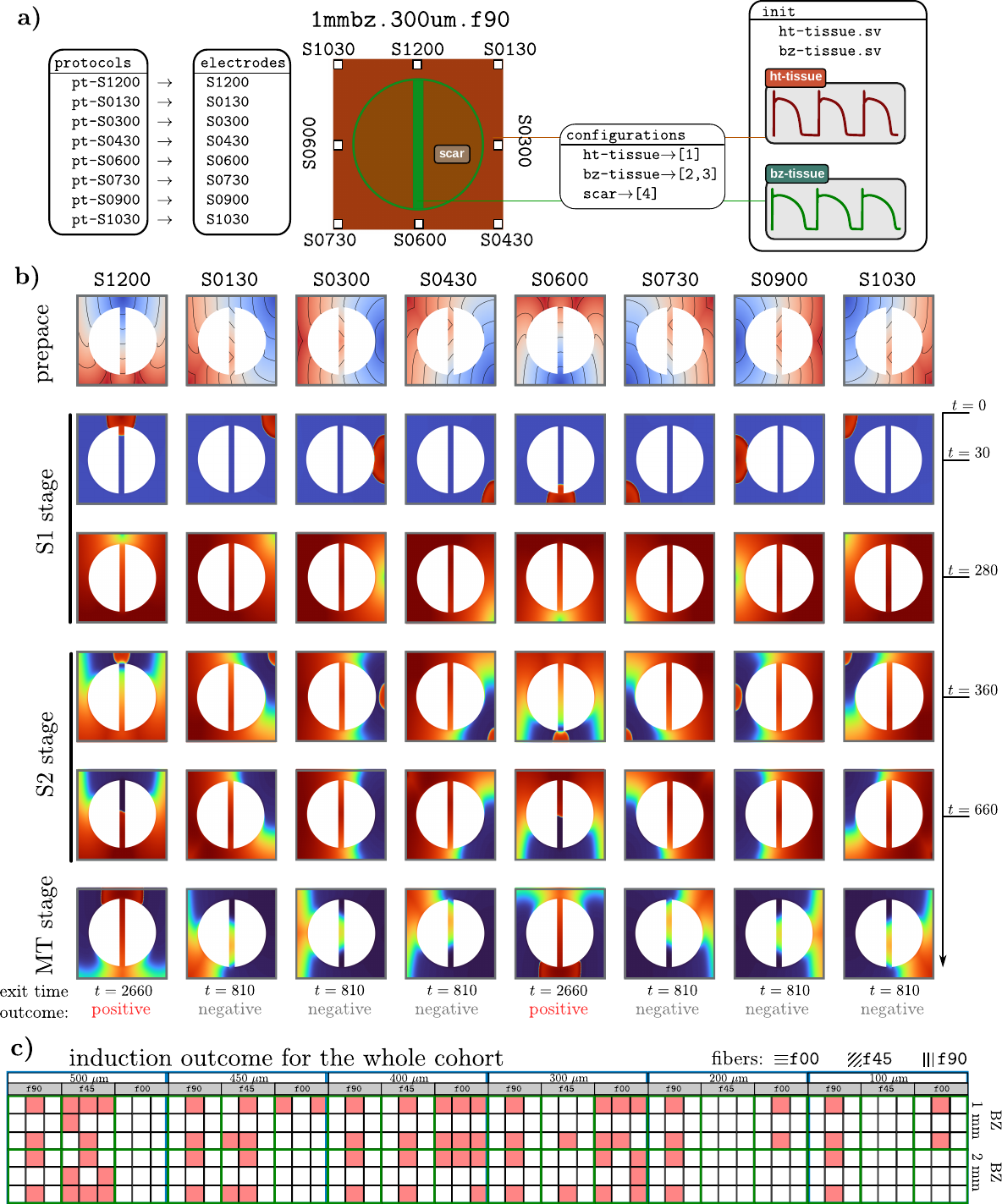}
    \caption{Example of experiment setup and results of the \aV pipeline applied subject \texttt{1mmbz.300um.f90} of the isthmus cohort, and results for the whole cohort. \textbf{a)} setup: eight protocols were defined around the isthmus and used to pace the tissue. Tissue properties were initialized using \oC and assigned according to the configurations section of the planfile. \textbf{b)} execution of \aV for subject  \texttt{1mmbz.300um.f90}. \gls{lat} maps computed during the prepace stage are shown, followed by two representative snapshots of S1 and S2 stages. The MT snapshots show the time when the simulation was interrupted, with the associated time stamp used to determine the outcome. \textbf{c)} table summarizing the outcome of \aV for the whole cohort considered. Red squares indicate electrodes which resulted in sustained reentry}  
    \label{fig:results_fig}
\end{figure*}

\noindent The whole process of \aV for subject \verb|1mmbz.300um.f90| of the cohort is shown in Fig. \ref{fig:results_fig}b.
For the selected subject, in the first row we display the \gls{lat} maps computed in the \gls{pp} stage, followed by the membrane potential saved in the checkpoints of the later \aV stages.
The images were generated by using the \verb+--movies+ option, or loading the \gls{lat} map of the \gls{pp} stage into Meshalyzer. Supplementary movies 1-8 show an animation of the transmembrane potential through the full induction protocol.  
The outcome of \aV for the whole cohort is shown in Fig. \ref{fig:results_fig}c, summarizing the results obtained for the whole cohort of subjects. Red squares highlight stimulation sites that resulted in sustained reentry for each subject. The orientation of the cells and the electrodes are the same as for the displayed subject.
For each subject, induction of reentry on each electrode was determined by reading the status tables generated with the \verb|--status| option, as shown shown in Tab. \ref{tab:status}.

To share the results of a specific stage of the experiment, it is sufficient to prompt again the same command used to generate the stage, adding the \verb|--bundle| flag. 
This will generate a zipped file containing all the files necessary to replicate the study, in the correct structure to perform \aV.  
The zipped file is a self-contained bundle with metadata and a default license file that is uploaded in the research data repository \href{https://radar.kit.edu/radar/de/home}{RADAR}, in compliance with the FAIR4RS principles for publication of research software \cite{houillon_facilitating_2024}. 
We provide the output of the \verb+--bundle+ option for subject \verb|1mmbz.300um.f90| as a Zenodo \href{https://zenodo.org/records/17159454}{record}.

\subsection{Test with anatomical geometry}
\label{sec:biv_setup}
\begin{figure*}[h!]
    \centering
    \includegraphics[width=.9\textwidth]{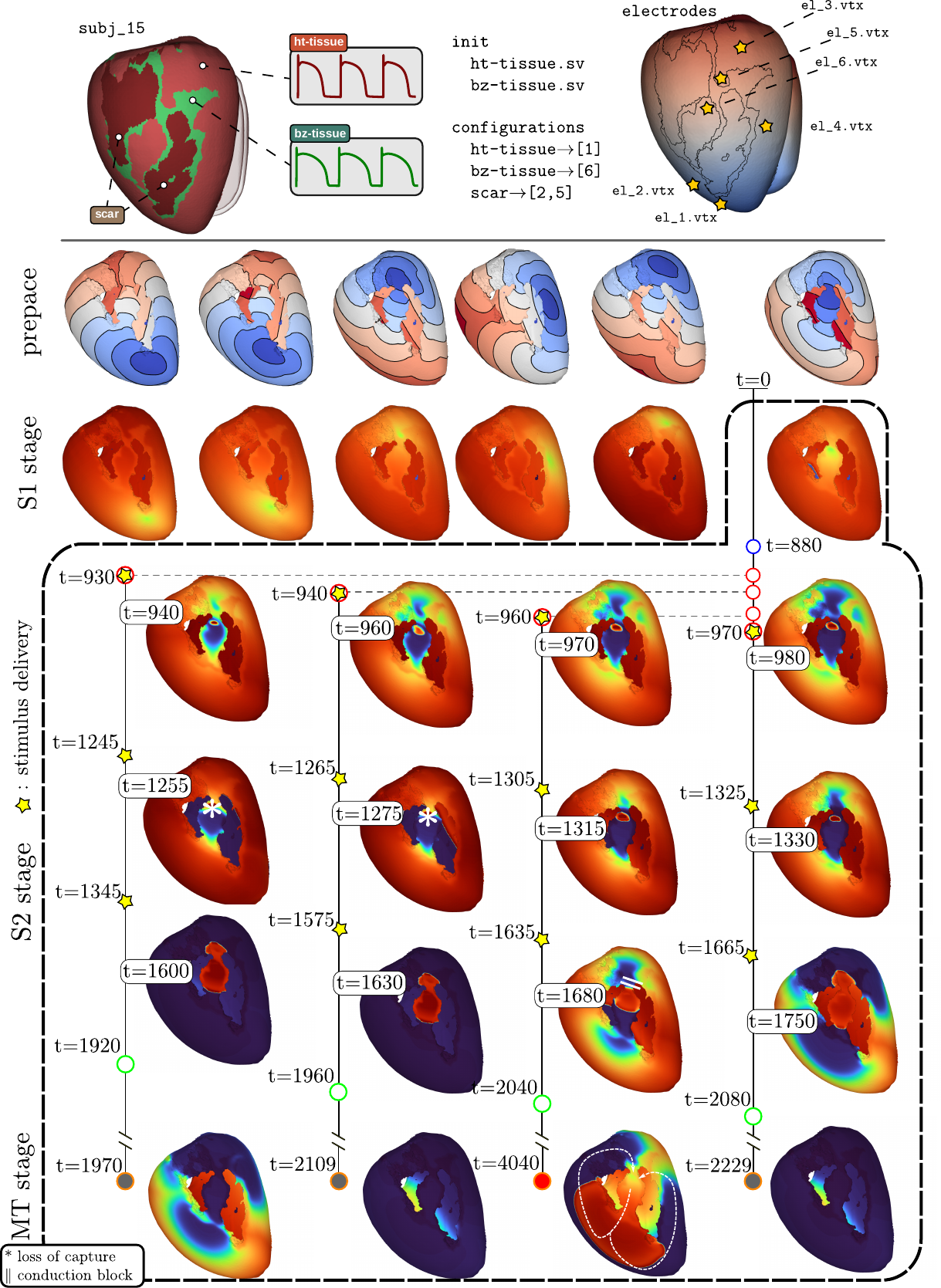}
    \caption{example of an \aV execution on an anatomical geometry. Top panel: experimental setup used to define the induction protocol.  
    Central panel: results of \aV for the \gls{pp} and S1 stage, computed for every electrode of the setup.
    We highlight the induction protocol for electrode 6 showing four different S2 intervals (330,340,360,370) delivered three times with a decrement of \SI{15}{\milli \second}. Early stimuli (330,340) fail to capture, whereas late stimuli elicit normal conduction in both healthy tissue and border zone. Unidirectional conduction block was achieved for a S2 interval of \SI{360}{\milli \second} and reveals the reentrant circuit, as assessed checking the exit time of the \gls{mt} stage. }
    \label{fig:results_fig_2}
\end{figure*}
To demonstrate the flexibility of \aV we briefly showcase an additional induction study using an anatomical geometry adapted from subject 15 in the open repository published in the work of Strocchi et al. \cite{strocchi_publicly_2020}. 
The original meshes are available in a resolution that is far too coarse to perform reaction-diffusion simulations, thus we resampled the geometry to  \SI{400}{\micro \meter} using \href{https://bitbucket.org/aneic/meshtool/src/master/}{meshtool}, and in house python code to interpolate \glspl{uvc}. 
After resampling, we performed tuning of tissue conductivities to compensate for numerical artifacts using the \verb|--tissue-tuning| option. 
The setup of the experiment is the same as for the previous case, with the exception that now we provide different electrodes and configurations files to adapt the workflow to our case. 
To save computational resources, the \gls{rv} was tagged as scar tissue in the configurations file so that it is not considered in simulations. 
Indeed, the presence of the \gls{rv} has an effect on \gls{vt} maintenance but the induction is local, so a test induction study can still be performed without it.
We also changed the number of electrodes used, so we provide a different planfile with a modified protocols section. In this example, we prescribe 5 electrodes using \glspl{uvc} coordinates and one electrode using mesh vertices to ensure that an induction site is paced. The tissue configurations and electrodes are shown in the top panel of Fig.  \ref{fig:results_fig_2}.

Since this setup required considerable computational resources, \aV was used to generate parameter files that were then uploaded to the \gls{hpc} facility Austrian Scientific Computing and used to launch simulations. To generate parameter files without running the experiment, we used the \verb|--gen-param-files| option. This allowed for an easy set up of the whole induction study just by uploading the parameter files in dedicated directories with the subject geometry. 
To reduce the use of computational resources, we only computed the \gls{mt} stage for the case of electrode 6, which was user-defined to ensure onset of \gls{vt}. Note that this does not affect the procedure in any way, since the parameters files relative to other electrodes were generated with the reported commands but were simply not simulated.\\
In this case we used the default prepacing options (\verb+--gen-lat ek --slim-cyc lat-1+) and decided to use two S1 stimuli to reach steady state since the geometry is more complex. 
Furthermore, we sweep four extrastimulus intervals of  330,340,360,\SI{370}{\milli \second} and use an aggressive pacing protocol of 3 S2 with a decremental interval of \SI{15}{\milli \second} provided through the option \verb|--decrement-S2 15|.
To generate all the files needed for the induction study, we executed \aV with the following options:\vspace{3mm}\\

\small
\begin{minipage}{\textwidth}    
\verb+python auto-varp.py --plan lf-vec-plan.json --cohort tmp-cohort --stage MT \+\\ 
\verb+    --configurations standard_VARP_configurations.json \+\\
\verb+    --electrodes standard_VARP_electrodes.json \+\\
\verb+    --gen-param-files --CI-array 330,340,360,370 --S1-cycles 2 --S2-cycles 3 \+\\
\verb+    --tissue_tuning --decrement 15+
\end{minipage}
\normalsize

\bigskip

The results of \aV applied to this case can be seen in Fig. \ref{fig:results_fig_2}. We show the setup and outcome of \gls{pp} and S1 stage for each electrode, then we focus on the induction study performed for electrode 6. For this electrode, we provide also the supplementary movies 9-12 showing the transmembrane potential through the whole induction protocol. 
Among the four tested intervals, only S2=\SI{360}{\milli \second} caused reentrant activity due to an unidirectional block at the entrance into the infarct isthmus.  

\section{Discussion}
Testing inducibility in virtual heart models has been used 
initially for investigating arrhythmia mechanisms,
with a more recent expansion of the scope towards clinical applications
such as non-invasive risk stratification of \gls{vt} recurrence\cite{arevalo16:_varp}, 
or guiding of \gls{vt} ablation therapy\cite{prakosa18:_vaat}.
The execution of such virtual inducibility studies is a complex and highly error-prone endeavor. 
Previous studies used proprietary in-house solutions, 
manually sweeping induction sites and pacing protocol parameters, 
involving interactive trial-and-error procedures
\cite{arevalo16:_varp,padilla_silico_nodate,camposAutomatedNearrealTime2022}.
This has been limiting accessibility of such studies to a small number of specialized labs
of highly trained users with expertise in computational modeling and simulation engineering.
This impeded a thorough cross validation of reported results by replication studies, 
typically considered an integral part of the scientific process.
However, this lack of reproducibility undermines the trustworthiness and credibility of study results 
that cannot be independently scrutinized, corroborated, or falsified outside the originating laboratory.
Moreover, design decisions in these studies aimed at reducing the high computational costs, 
favoring speed over accuracy to keep studies tractable. 
These design decisions have been shown to be associated with a risk of introducing stochasticity, 
raising concern about potentially spurious results \cite{bishop25:_stochastic}. 
For further advancing virtual inducibility testing to comply with industrial and clinical safety requirements \cite{bodnerFrameworkSilicoClinical2022}
a broader community effort on validating this approach appears warranted, 
to better delineate the bounds of applicability of this powerful methodology.

Here, we address this demand by reporting on the implementation of a computational pipeline, \aV, 
that expands accessibility of this technology to researchers with electrophysiological expertise, 
without the need for years of training in simulation engineering.
The \aV pipeline and all its integral components required for execution, analysis and post-processing
are built upon openly accessible tools provided through the \oC project \cite{plank21:_opencarp}, 
with \oC, \cU and \fC being the key components.

\subsection{Efficient scalable inducibility testing}

Achieving efficiency and scalability of workflows must address the two major cost factors 
in virtual inducibility studies. These are the computational costs of execution the study and, 
most importantly, 
the human operator costs incurring by defining, coding, launching, monitoring and analyzing 
the complex processing workflows.
An important factor is a high level of abstraction
that facilitates the application of the protocol to any type of anatomical model,
using an arbitrary number of electrode positions for testing. 
In \aV, all experimental and control parameters are compactly exposed as a set of dictionaries in \verb+.json+ format.
These can be edited to modify model parameters, electrode locations or protocol definitions, 
globally or at the individual subject level, thus separating experimental definition from execution.
This enforces a standardized consistent execution of the protocol over an entire cohort,
ascertaining that inducibility in each virtual subject is tested in the same way, 
to avoid spurious test outcomes as they inevitable occur with manual interrogation.
Importantly, the translation of abstract protocol definitions to executable simulation commands 
is steered by providing high level arguments,
thus reducing errors by alleviating users from any coding tasks. 

The execution of the \aV protocol is highly automated. 
In principle, \aV can be applied to a single subject or an entire virtual cohort, 
by a single command without prompting for any further user input.
The generation of all commands and parameter files 
required for performing an entire induction study readily supports the execution 
using \gls{hpc} facilities, if required by the workload.
 
Further, the computational efficiency of the \aV pipeline has been optimized 
by using most recent advanced features of the underlying \oC environment. 
Substantial computational savings are achieved through a more efficient calculation of the reference limit cycle 
using \fC \cite{gsell24:_forcepss} as compared to previously employed brute-force approaches\cite{arevalo16:_varp}. 
For instance, these studies simulated 7 seconds of S1 pacing to arrive at a stable limit cycle,
with simulations lasting for hours when using hundreds of compute cores 
at high performance computing facilities. 
In our implementation this is achieved in $\simeq$1.2 seconds on a desktop computer
by using the \fC limit cycle feature\cite{gsell24:_forcepss}.
Further savings are gained from a staged execution
that separates the pre-conditioning phases in stages \gls{pp} and S1 
from the induction and maintenance testing phases in stages S2 and \gls{mt}.
This ensures that the computationally most expensive simulation phases are reused wherever possible,
and that costly repetition of phases overlapping between variations of the protocol are avoided.
The most important gain in efficiency though is due to the reduced need for operator input,
as tedious progress monitoring and interactive evaluation of intermediate results is avoided.

\subsection{Reproducibility and validation}
Previously reported virtual inducibility studies using the VARP protocol \cite{arevalo16:_varp,prakosa18:_vaat,sung22:_inFAT} are, 
without exception, non-reproducible for a number of reasons:
i) The entire complex processing workflow and its implementation is not exposed or concisely described;
ii) The numerical solver scheme and the bulk of numerical parameters influencing simulated behavior 
have not been reported. In a complex simulator such as \oC many schemes are implemented (implicit/explicit, operator splitting, mass lumping, stiffness assembly, etc) 
and a broad range of parameters selecting solver and convergence metrics can be used. 
Any choice may influence the outcome of the inducibility experiment, 
particularly when coarser mesh resolutions are used for the sake of saving computational cost\cite{bishop25:_stochastic}.
iii) Computational grids, pacing locations and simulation outputs are not made available.
Thus, an evaluation of study outcomes is effectively impeded. 
As such, the validity of virtual inducibility studies can be believed or not, 
but not verified or falsified which is unsatisfying from a scientific point of view.

These limitations are addressed by our implementation of \aV, 
with the exception of model and data accessibility, 
which remains at the user’s discretion.
The \aV workflow generates a comprehensive audit trail which uniquely documents 
all model parameters and how these are processed to produce the study results.
All integral software components used in simulation, data analysis and visualization 
are freely available and accessible.
This renders studies reproducible by using \aV for generating simulation data, 
without the need of reverse engineering processes and parameters that are not fully disclosed in publications,
given that all, or at least selected subsets of computational meshes are made available as well.
To minimize the additional costs a tool for uploading study data in a FAIR repository is integrated
using the \textit{bundle} option of \oC.
Thus, if computational meshes are made available studies using \aV are identically reproducible by others,
offering a basis for cross validation of studies within the cardiac modeling community.


\subsection{Robustness}
The \aV implementation enhances robustness by ensuring a correct execution of the study, 
by avoiding user errors inevitably occurring with manual steering of complex simulation pipelines.
Indeed, it is important to ensure robustness of predictions with respect to uncertainties in input parameters.
Simulations where a minor perturbation of an uncertain input parameter entails a marked change in prediction
are not robust and should not be trusted.
An inducibility test such as \aV is only robust 
if predictions on inducibility, maintenance and behavior of induced arrhythmias 
are insensitive to input parameter uncertainty.
As shown previously \cite{mendoncacostaDeterminingAnatomicalElectrophysiological2022,bishop25:_stochastic}, 
this is not the case as the outcome of induction studies is highly sensitive to input uncertainty and modeling choices. 
For instance, altering the width of the border zone around a scar by 1 mm 
may yield vastly different outcomes (refer to Fig. \ref{fig:results_fig}). 
Also technical factors such as spatial mesh resolution play a major role 
by influencing the speed of wave front propagation.
In particular, this is the case for slow conduction associated with steep wave fronts in space 
which are challenging to adequately resolve.
Thus, the outcome of testing inducibility at mesh resolution of \SI{500}{ \micro \meter} 
may be entirely different from the same test at a resolution of \SI{100}{ \micro \meter} (refer to Fig.\ref{fig:results_fig}).

With \aV this issue is not directly addressed, 
but abstraction and scalability of the \aV framework supports their investigation. 
More robust testing can be achieved with \aV by computing a set of solution per individual model 
where each solution is obtained by perturbing the input parameters.
This may help to ascertain that a prediction robustly depends on the model's anatomy and physiology,
and not on random uncertainties.
This is exemplified in an example cohort of geometrically simple infarct models comprising isthmus and outer loop
(Fig.~\ref{fig:results_fig}).
When modeling an infarcted patient heart the exact geometry of the scar and border zone 
are determined by segmentation which is afflicted with uncertainty. 
Fiber orientations cannot be measured accurately \emph{in vivo} 
and are approximated by using rules empirically based on \emph{ex vivo} observations \cite{rohmer07:_reconstruction,bayer12:_rule_based}.
Keeping anatomy of infarct and tissue constant, we tested the impact of 
the width of the border zone by only 1 mm which is within the segmentation uncertainty,
the uncertain fiber orientation with respect to the orientation of the isthmus, 
as well as spatial resolution.
Accounting for this uncertainty revealed that the outcome of the test is entirely random, 
and highly sensitive to input parameters.
For instance, at the highest resolution grid with a mesh resolution of \SI{100}{ \micro \meter} 
the model is inducible for fibers aligned with the isthmus, \verb|f90|, for a border zone width of \SI{1}{ \milli \meter}. 
With a minimal change in border zone width to \SI{2}{\milli \meter} or altering fiber angle to \verb|f45| 
the model becomes non-inducible. 
Consistency between spatial resolutions is only witnessed between the \SI{100}{ \micro \meter} and 
\SI{200}{ \micro \meter} grids whereas inconsistent and contradictory results are seen 
at all higher resolutions. 

While one may argue that a tissue calibration step like the one performed in section \ref{sec:biv_setup} could solve these issues, restitution properties of propagation velocity may still cause simulated velocities to be well below the acceptable threshold for the given resolution, resulting in unpredictable behavior.
On the other hand, mesh resolutions are constrained by computational resources allocated to a specific study and often a lower bound for which studies become completely unfeasible is found.
Determining the bounds within numerical schemes are appropriate for a specific study is well outside the scope of our work, however we partially address this matter by providing all the tools necessary to foster sharing and scrutiny of studies between groups. 

\subsection{Limitations}
For the sake of computational efficiency the current implementation of \aV supports only a standard S1-S2 induction protocol.
However, inducibility also depends on details of the induction protocol employed. 
Using a more sophisticated protocol with arbitrary number of stimuli and coupling intervals 
may achieve induction where a simpler S1-S2 protocol fails to induce.
As such, a more fine-grained testing of protocols may therefore alter the prediction on inducibility.
However, as in clinical settings, an exhaustive parameter sweep may substantially increase costs, making it challenging to keep the tests computationally tractable.
Using more elaborate protocols with \aV is feasible, but requires additional coding.

Further, support for quantifying sensitivity of predictions to input uncertainty has not been implemented in \aV.
Considering the vast dimensionality of the model itself, with a large number of highly uncertain parameters
that govern electrophysiological behaviors at the cellular and tissue scale,
the outcome of virtual inducibility testing with \aV must be considered probabilistic.
However, this has not been explored in any of the virtual inducibility studies reported in the literature.
A thorough quantification of how parameter and initial condition uncertainty impacts inducibility outcomes 
is essential for providing estimates of confidence. 
The scalability of \aV supports such analysis by deriving perturbed models for each case in a cohort,
as demonstrated for the geometrically simplified isthmus cohort (Fig.~\ref{fig:results_fig}), 
but manual parameter sampling and output processing is required.

\FloatBarrier

\section{Conclusion}
The success of emerging industrial and clinical applications of computational modeling 
such as testing of arrhythmia inducibility in virtual hearts
relies upon implementation of \emph{in silico} experiments that are robust and reproducible, 
and can be subjected to scrutiny and a rigorous assessment of validity.
This gap is addressed by providing \aV as a unified framework for virtual induction studies,
and showcasing its application for studying scar-mediated \glspl{vt}.
The \aV workflow achieves a high degree of automated processing, requires minimal user input 
and avoids explicit coding which standardizes the procedure, reduces operator errors 
and expands the user base beyond expert simulation engineers.
This renders inducibility testing scalable as test are readily executed across virtual cohorts of subjects
in an unattended manner.
Furthermore, \aV automatically generates an audit trail used to track model parameters and processing
steps that can be used to reproduce study results. This is further facilitated 
by the possibility of uploading study data in a FAIR repository thanks to the 
\textit{bundle} option of \oC. By providing computational meshes used in the study,
computational experiments become completely reproducible.



\newpage
\appendix
\section{Structure of \gls{uvc} based electrodes and template plan file}
\label{sec:appendix:plan}

Here we provide a brief example of a simulation plan file and show the complete definition of \gls{uvc} based electrodes, which were not available in the initial \fC implementation. For an exhaustive definition of a plan file refer to the original \fC publication \cite{gsell24:_forcepss}, or to the \href{https://zenodo.org/records/17159454}{online supplementary material} \\

\noindent
\verb|"functions":{|\\
\verb|  "version": 2,|\\
\verb|  "definitions":{|\\
\verb+    "ht_tissue":{    +\\
\verb+      "ht_tissue":{+\\
\verb+        "model":                "tenTusscherPanfilov"+\\
\verb+        "model_par":            null +\\
\verb+        "plugins":              null+\\
\verb+        "plugins_par":          null+\\
\verb+        "initialization":{    +\\
\verb+          "num_cycles":         100+\\
\verb+          "bcl":                600+\\
\verb+          "init":               null+\\
\verb+          "apdres_file":        null+\\
\verb+          "apdres_protocol":    null+\\
\verb+        }    +\\
\verb+      },   +\\
\verb+      "conductivity":{    +\\
\verb+        "gil":    0.255+\\
\verb+        "gel":    0.625+\\
\verb+        "git":    0.0775+\\
\verb+        "get":    0.236+\\
\verb+        "gin":    0.0775+\\
\verb+        "gen":    0.236+\\
\verb+        "surf2vol":    0.14+\\
\verb+      },    +\\
\verb+      "conduction_velocity":{    +\\
\verb+        "reference":{    +\\
\verb+          "vf":    0.6+\\
\verb+          "vs":    0.2+\\
\verb+          "vn":    0.2+\\
\verb+        },    +\\
\verb+        "measured":{    +\\
\verb+          "vf":    null+\\
\verb+          "vs":    null+\\
\verb+          "vn":    null+\\
\verb+        }    +\\
\verb|      }|\\
\verb|    },|\\
\verb+    "bz_tissue":{    +\\
{\color{red}\verb+     ...   +\\}
\verb+    }    +\\
\verb+    "scar":{    +\\
{\color{red}\verb+     ...   +\\}
\verb+    }    +\\
\verb+  }    +\\
\verb+},    +\\
\verb|"protocols":{|\\
\verb|  "version":2,|\\
\verb|  "prepacing":{|\\
\verb|    "protocol_1":{|\\
\verb|      "propagation": "rd"|\\
\verb|      "num_cycles": 1|\\
\verb|      "bcl": 600|\\
\verb|      "electrodes":|{\color{green}\verb|"uvc-el"|}\\
\verb|      "rel_timings": null|\\
\verb|      "lat_file": null|\\
\verb|      "restart": null|\\
\verb|    }  |\\
\verb|  } |\\
\verb|},|\\
\verb|"electrodes": {|\\
\verb|  "version": 2,|\\
\verb|  "definitions": {|\\
\verb+    "+{\color{green} \verb+uvc-el+}\verb+": +{\color{green}\verb+{+}\\
\verb|      "type": "ucc_sphere",|\\
\verb|      "p0": list(type=float, length=3, desc="uvc point, apicobasal, transmural, and rotational")|\\
\verb|      "cavity": value(type=str,choices={"lv","rv"}),|\\
\verb|      "radius": value(type=float, unit=mm, desc="threshold value defined in meshtool"),|\\
\verb+      "searchdom": "cxyz"+\\
{\color{green}\verb|    }|}\verb+,+\\
{\color{red}\verb|    ...|}\\
\verb+  }+\\
\verb+}+\\
\verb|"configurations": {|\\
\verb|  "version":2,|\\
\verb|  "definitions":{|\\
\verb|    "healthy_tissue":{|\\
\verb|      "tags": [1],|\\
\verb|      "func":"ht_tissue"|\\
\verb|    },|\\
\verb|    "borderzone_tissue": {|\\
{\color{red}\verb|      ...|\\}
\verb|    },|\\
\verb|    "scar": {|\\
{\color{red}\verb|      ...|\\}
\verb|    }|\\
\verb|  }|\\
\verb|},|\\
\verb|"solver_setup": {|\\
{\color{red}\verb|  ...|\\}
\verb|}|\\
Green keywords are exemplary and are user-defined. Red dots indicate missing sections that must be filled according to \fC standards. 
In the \verb+electrodes+ section, the values of individual entries are specified as predefined keywords (e.g., \verb|"lv"|), single values, or lists of values. Single values are declared by the \verb|value| keyword and specified by their type, unit, and short description. Lists of values are declared by the \verb|list| keyword and specified by their length, type of entries, and a short description. More details can be found in the original \fC work \cite{gsell24:_forcepss}

\section{Naming convention of checkpoint files}
\label{sec:appendix:checkpoints}
Here we provide the naming conventions for the checkpoint files saved in each stage of \aV.  
Each checkpoint files refers to user defined inputs that are described here:\\
\textit{protocol}: name of the protocol definition in the plan file or protocols file.\\
\textit{gen-lat}: option used to generate the activation times file, either 'rd' or 'ek'.\\
\textit{electrode}: name of the electrode definition used by the protocol.\\
\textit{meshname}: name of the folder containing the subject geometry and accessory files.\\
\textit{lim-cyc}: option provided to compute the limit cycle stage in the \gls{pp} stage.\\
\textit{pcl}: pacing cycle length used when computing the steady state of the experiment.\\
\textit{exit\_time}: time at which the simulation is stopped.\\
\vspace{4mm}
\textit{S2\_interval}: extrastimulus interval considered for the simulation (only used in S2 stage simulations).\\
\smallskip
Naming convention of checkpoint files, separated by stage:\\
lat-based prepacing activation sequence:\\
\smallskip
\verb|lim_cic-ptcl_{}_{}_el_{}-act_seq.dat| (\textit{protocol,gen-lat,electrode})\\
prepacing checkpoint file:\\
\smallskip
\verb+{}_{}_pp_{}_bcl_{:.1f}_tstamp_{:.1f}.roe+ (\textit{protocol, meshname, lim-cyc, pcl, checkpoint\_time})\\
S1 checkpoint file:\\
\smallskip
\verb|S1_{}_PCL_{:.1f}_ms_{}_tstamp_{:.1f}.roe| (\textit{electrode,pcl,meshname, exit\_time})\\
S2 checkpoint file:\\
\smallskip
\verb|S2_{}_PCL_{:.1f}_ms_CI_{:.1f}_ms_{}_tstamp_{:.1f}.roe| (\textit{electrode,pcl,S2\_interval,meshname,exit\_time})\\
MT checkpoint file:\\
\verb+MT_{}_PCL_{:.1f}_ms_CI_{:.1f}_ms_{}_tstamp_{:.3f}.roe+ (\textit{electrode,pcl,S2\_interval,meshname,exit\_time})

\printcredits

\bibliographystyle{unsrt}

\bibliography{auto-varp}
\end{document}